\documentclass[10pt, twoside, leqno]{elsarticle}
\usepackage{amsmath, amssymb, amsthm, epsfig}
\usepackage{hyperref, latexsym}
\usepackage{url}
\usepackage[mathscr]{euscript}
\usepackage{pdfpages}
\usepackage{color}
\usepackage{fullpage}
\usepackage{setspace}

\usepackage[enableskew]{youngtab}
\def\KeyWord#1{$\backslash$\IfColor{$\!\!$\textRed{#1}\textBlack}{#1}$\!\!$}

\usepackage{amsfonts}
\usepackage{amssymb}
\usepackage{amsmath}
\usepackage{amsthm}

\usepackage{etex}
\usepackage{amsgen}
\usepackage{amsopn}
\usepackage{verbatim}
\usepackage{xspace}
\usepackage{multicol}
\usepackage{eepic}
\usepackage{upref}
\usepackage{pgf}
\usepackage{tikz}
\usetikzlibrary{patterns}
\usepackage{pgffor}
\usepackage{pgfcalendar}
\usepackage{pgfpages}
\usepackage{shuffle,yfonts}

\def\today{\ifcase\month\or
  January\or February\or March\or April\or May\or June\or
  July\or August\or September\or October\or November\or December\fi
  \space\number\day, \number\year}

 \theoremstyle{definition}

 \theoremstyle{remark}

\newtheorem{thm}{Theorem}
\numberwithin{equation}{section} 

\begin{document}

\title{All of zeros of Riemann's Zeta-Function are on $\sigma$=1/2}
\author{Nianrong Feng}
\ead{fengnr@mail.ahnu.edu.cn}
\address{School of Computer and Information, Anhui Normal University}
\author{Yongzheng Wang}
\address{School of Mathematics and Statistics, Anhui Normal University}

\numberwithin{equation}{section}
\begin{frontmatter}
\begin{abstract}
The research shows that Riemann proved that all of zeros of Riemann's zeta function are on $\sigma=1/2$ based on the functional equation
\begin{align*}
 \pi^{-\frac{s}{2}}\Gamma \left( \frac{s}{2}  \right) \zeta(s)&={\frac{1}{s(s-1)} + \int\limits_1^\infty \psi(x) \left(
            x^{\frac{s}{2} - 1} + x^{-\frac{1+s}{2}}
            \right) \,dx,}\quad\qquad{s}=\sigma+it,
\end{align*}
which is in Riemann's ``\"{U}ber die Anzahl der Primzahlen unter einer gegebenen Grosse". According to the geometric meaning of the functional equation and the argument principle, we obtain the number of zeros $N_0(T)$ of the Riemann zeta function on the critical segment $\sigma=1/2,0\leq{t}\leq{T}$ and the number of zeros $N(T)$ of the Riemann zeta function in the rectangular region $-1\leq\sigma\leq{2},0\leq{t}\leq{T}$, respectively. The result is
\begin{align*}
N(T)&=N_0(T)=\frac{\arg{\left[\pi^{-\frac{s}{2}}\Gamma\left(\frac{s}{2} \right)\zeta(s)\right]}}{\pi}+1\\
&=\frac{T}{2\pi}\log\frac{T}{2\pi}-\frac{T}{2\pi}+O(\log{T}),\qquad{s=1/2+iT}.
\end{align*}
\end{abstract}

\end{frontmatter}

\section{Notations}
Let $s= \sigma+it$ be a complex number, $\zeta(s)$ be the Riemann zeta-function.
$\mathbf{D(T)}$ and $\mathbf{L(T)}$ are defined as
\[\mathbf{D}(T)=\{\sigma+it:-1\leq\sigma\leq{2},0\leq{t}\leq{T},\zeta(\sigma+iT)\neq{0}\},\]
\[\mathbf{L}(T)= \{\sigma+it:\sigma= 1/2,0\leq{t}\leq{T},\zeta(1/2+iT)\neq{0}\}. \]
Let $N_{0}(T)$ be the number of zeros of $\zeta(1/2+it)$ on $\mathbf{L}(T)$, and $N(T)$ be the number of zeros of $\zeta(s)$ in $\mathbf{D}(T)$. The essence of the Riemann Hypothesis is to prove that $N(T)=N_0(T)$.
\section{Introduction}
The reference \cite{riemann} is Riemann's only paper on number theory, in which he studied some properties of the Riemann zeta function and the relationship between its zero point and the distribution of prime numbers, and for the first time studied number theory with the ideas and methods of Complex Variable Function. The reference \cite{riemann} consists of three paragraphs, each of which has a conclusive formula, which is in turn:
\begin{align}\label{Pxzgsx5x}
{\Gamma \left( \frac{s}{2} \right) \pi^{-\frac{s}{2}} \zeta(s)=\Gamma \left( \frac{1-s}{2} \right) \pi^{-\frac{1-s}{2}} \zeta(1-s)}
\end{align}
\begin{align}\label{sggsx1}
\pi^{-\frac{s}{2}}\Gamma(\frac{s}{2})\zeta(s)&=\frac{1}{s(s-1)}+\int_1^{\infty}\psi(x)(x^{\frac{s}{2}-1}+x^{-\frac{s+1}{2}})dx, \psi(x)=\sum_{n=1}^\infty{e^{-n^2\pi{x}}},
\end{align}
\begin{align}\label{CXGS}
F(x)&=\sum{(-1)}^\mu\frac{1}{m}J(\frac{1}{m}),
\end{align}
where
\begin{align*}
J(x)&=Li(x)-\textcolor[rgb]{1.00,0.00,0.00}{\underline{\sum^\alpha\left(Li(x^{1/2+\alpha{i}})+Li(x^{1/2-\alpha{i}})\right)}}+\int_x^\infty{\frac{1}{x^2-1}\frac{dx}{x\log(x)}}+\log\xi(0),
\end{align*}

Academics have explained the functions of these three formulas. Formula \eqref{Pxzgsx5x} is the functional equation for the Riemann zeta function. Formula \eqref{sggsx1} is the second proof of the functional equation for the Riemann zeta function\cite{Siegel,Titch,xref1}. Formula \eqref{CXGS} is a formula for calculating the number of prime numbers based on the Riemann Hypothesis. It should be emphasized that, so far, thousands of literatures have taken the Riemann Hypothesis as a condition to prove the problems related to the zero of Riemann zeta function\cite{ref19}.

According to the above explanation, this formula \eqref{CXGS} is a meaningless formula because it is obtained based on the greatest unsolved problem. According to common sense, Riemann will not go against basic common sense to do meaningless work. Our research shows that the source of this contradiction is a one-sided understanding of the formula \eqref{sggsx1}. Here, it is not the second proof of the functional equation of the Riemann zeta function, but the geometric proof of the proposition that \textbf{all nontrivial zeros of $\zeta(s)$ function lie on the line Re (s)=1/2}. In \cite{riemann}, Riemann made a technical treatment to this proof, that is, it shielded the most critical geometric formula of $\zeta(s)$ function. After adding this formula, we can clearly see that he proves this proposition from two aspects.

\begin{description}
  \item[The first aspect:] According to the formula \eqref{sggsx1}, $s=0,1$ are the two poles of $\pi^{-\frac{s}{2}}\Gamma(\frac{s}{2})\zeta(s)$. Set
$$\xi(s)=\frac{s(s-1)}{2}\pi^{-\frac{s}{2}}\Gamma(\frac{s}{2})\zeta(s).$$
\end{description}
According to the argument principle, the number of zeros of $\xi(s)$ function in the rectangular area $D(T)$ is approximately
\begin{align*}
N(T)&\approx\frac{T}{2\pi}\log\frac{T}{2\pi}-\frac{T}{2\pi}.
\end{align*}

\begin{description}
  \item[The second aspect:] By the formula \eqref{sggsx1}, when $s=1/2+it$,
$$\xi(s)=\frac{1}{2}-(t^2+\frac{1}{4})\int_{1}^{\infty}\psi(x)x^{-\frac{3}{4}}\cos(\frac{1}{2}t\log{x})dx,$$
\end{description}
and for any $s=1/2+it$, $\xi(s)$ function is a finite real function. \textcolor[rgb]{0.00,0.00,0.50}{\textbf{
This conclusion shows that $\pi^{-\frac{s}{2}}\Gamma(\frac{s}{2})\zeta(s)$ function must have the geometric formula}
\begin{align}\label{jsggsx1}
\mathbf{\pi^{-\frac{s}{2}}\Gamma(\frac{s}{2})\zeta(s)}&=\mathbf{\varphi(s)+\varphi(1-s),}
\end{align}
\textbf{and when $s=1/2+it$}
\begin{align}\label{jsggsx12}
\mathbf{\pi^{-\frac{s}{2}}\Gamma(\frac{s}{2})\zeta(s)}&=\mathbf{2\Re\left[\varphi(s)\right]}\nonumber\\
&=\mathbf{2\left|\varphi(s)\right|\cos\left[\arg\varphi(s)\right].}
\end{align}
\textbf{When $t={T}$, according to $\cos(k\pi-\frac{\pi}{2})=0$ and
$\left|\arg\left[\pi^{-\frac{s}{2}}\Gamma(\frac{s}{2})\zeta(s)\right]-\arg\varphi(s)\right|<\frac{\pi}{2}$ , It is obtained from this formula \eqref{jsggsx12} that the number of zeros of $\zeta(s)$ function on the critical line $L(T)$ is equal to }
\begin{align}\label{xccsc}
N_0(T)=N(T)-C,\qquad{C}\textbf{ be a constant.}
\end{align}}
This result is consistent with Riemann's statement that \textbf{One now finds \textcolor[rgb]{1.00,0.00,0.00}{indeed} approximately this number of real roots within these limits, and it is very probable that all roots are real}\cite{riemann}.

The main reason for Riemann shielding this formula \eqref{jsggsx1} is that he encountered difficulties in proving its existence.
C. L. Siegel\cite{Siegel} said that in the first draft of the reference \cite{riemann}, there was a formula
\begin{align*}
\pi^{-\frac{s}{2}}\Gamma(\frac{s}{2})\zeta(s)&=\varphi(s)+\varphi(1-s)+R(s),
\end{align*}
where $R(s)$ is a remainder term ,and $$\varphi(s)=\pi^{-\frac{s}{2}}\Gamma(\frac{s}{2})\sum_{n=1}^{M}{\frac{1}{n^{s}}}\qquad{(s=\sigma+it,M=\lfloor\sqrt{\frac{t}{2\pi}}\rfloor)}.$$
Obviously, it has a remainder $R(s)$ that is difficult to eliminate, which is an approximate formula of the formula \eqref{jsggsx1}. In reference \cite{riemann}, Riemann stated the reason for deleting this formula, that is, `` Certainly one would wish for a stricter proof here; I have meanwhile temporarily put aside the search for this after some fleeting$\cdots$"\footnote{Unfortunately, the academia misunderstood this sentence and took it as a strong evidence for Riemann to propose the Riemann Hypothesis.}.
Fortunately, Riemann finally solved the problem of the existence of this formula \eqref{jsggsx1}.

Based on the above analysis, we summarize Riemann's brief proof of $N(T)=N_0(T)$ into four theorems, and make supplementary proofs for two of the key theorems.

\section{Four Theorems to Prove Riemann Hypothesis}
\begin{thm}\label{xdl1}
Assume $\Phi(s)=\pi^{-s/2}\Gamma(\frac{s}{2})\zeta(s)$, $\xi(s)=\frac{s}{2}(s-1)\Phi(s)$. There must be
\begin{align}\label{jhdcx1}
\Phi(s)&=\Phi(1-s),
\end{align}
and $\xi(s)$ is an entire function.
\end{thm}
The theorem proves that the zeros of $\xi(s)$ function in the strip $0\leq\Re(s)\leq{1}$ are symmetric about the vertical line $\Re(s)=1/2$, and also proves that the zeros of $\xi(s)$ function are precisely the non-trivial zeros of $\zeta(s)$ function.

This is the key theorem to prove that $N(T)=N_0(T)$. In \cite{riemann}, Riemann proves the theorem in detail, which will not be repeated in this paper.

\begin{thm}\label{xdl2}
The number of zeros of $\xi(s)$ function in the rectangular area $\mathbf{D}(T)$ is equal to
\begin{align}\label{NTF1}
\textcolor[rgb]{1.00,0.00,0.00}{\mathbf{N(T)=\frac{\theta(T)+\arg{\zeta(\frac{1}{2}+T)}}{\pi}+1}}.
\end{align}
\end{thm}
The expansion of formula \eqref{NTF1} is Riemann-von Mangoldt formula
\begin{align}\label{VONFORM}
{N(T)=\frac{T}{2\pi}\log\frac{T}{2\pi}-\frac{T}{2\pi}+O(\log{T})}.
\end{align}
The reason why the formula \eqref{VONFORM} is not used to represent $N(T)$ in this paper is that $O(\log{T})$ cannot perform the four fundamental operations of arithmetic.
In view of the fact that $N(T)$ given in \cite{riemann}  is an approximate value, which does not meet the requirements of accurate operation, this paper will
reprove the theorem by using the argument principle, and give the exact value of $N(T)$, that is, the formula \eqref{NTF1}.
\begin{thm}\label{xdl3}
It's known that the Riemann function equation
\begin{align}\label{xgs1}
 \pi^{-\frac{s}{2}}\Gamma \left( \frac{s}{2}  \right) \zeta(s)={\frac{1}{s(s-1)} + \int\limits_1^\infty \psi(x) \left(
            x^{\frac{s}{2} - 1} + x^{-\frac{1+s}{2}}
            \right) \,dx.}
\end{align}
From its geometric meaning, it is obtained that the number of non-trivial zeros of $\zeta(s)$ function on the critical line $\mathbf{L(T)}$ is
\begin{align*}
\textcolor[rgb]{1.00,0.00,0.00}{\mathbf{N_0(T)=\frac{\theta(T)+\arg{\zeta(\frac{1}{2}+iT)}}{\pi}+1}}.
\end{align*}
\end{thm}
\textbf{The geometric meaning of this formula \eqref{xgs1} is exactly formula \eqref{jsggsx1} that Riemann hoped to get}. Because
\begin{align*}
\pi^{-\frac{s}{2}} \Gamma \left( \frac{s}{2}  \right) \zeta(s)&=-\frac{1}{s} -{\frac{1}{1-s}+ \int\limits_1^\infty \psi(x) x^{\frac{s}{2} - 1} \,dx+ \int\limits_1^\infty \psi(x) x^{\frac{1-s}{2}-1} \,dx}\\
&=\left[-\frac{1}{s} + \int\limits_1^\infty \psi(x) x^{\frac{s}{2} - 1} \,dx\right]+\left[-{\frac{1}{1-s}+ \int\limits_1^\infty \psi(x) x^{\frac{1-s}{2}-1} \,dx}\right]
\end{align*}
Let
$$\varphi(s)=-\frac{1}{s}+\int_1^{\infty}\psi(x)x^{\frac{s}{2}-1}dx,$$
then
\begin{align*}
\pi^{-\frac{s}{2}}\Gamma(\frac{s}{2})\zeta(s)=\varphi(s)+\varphi(1-s).
\end{align*}

If Riemann had noticed the method of transformation, we would have guessed that there would be neither the following theorem nor the Riemann Hypothesis.
\begin{thm}\label{xdl4}
$\zeta(s)$ function has a functional equation in the form of
\begin{align}\label{xgs2}
\pi^{-\frac{s}{2}}\Gamma\left(\frac{s}{2}\right)\zeta(s)=\varphi(s)+\varphi(1-s),
\end{align}
and
\begin{align}\label{xgs21}
\pi^{-\frac{s}{2}}\Gamma(\frac{s}{2})\zeta(s)=2\Re\left[\varphi(s)\right],\qquad\Re(s)=1/2,
\end{align}
where
$$\varphi(s)=\pi^{-\frac{s}{2}}\Gamma\left(\frac{s}{2}\right)f(s), \qquad
f(s)=\int_{0\swarrow{1}}\frac{x^{-s}e^{\pi{ix^2}}}{e^{\pi{ix}}-e^{-\pi{ix}}}dx.$$
\end{thm}
This is a theorem proving the existence of this formula \eqref{jsggsx1}, where the formula \eqref{xgs2} is the famous Riemann-Siegel integral formula. In reference \cite{Siegel}, Riemann purposely emphasized that when $s=1/2+it$, the formula \eqref{xgs21} is obtained from the formula \eqref{xgs2}, which is exactly the same as the statement that we emphasized earlier that this formula \eqref{jsggsx12} is obtained from this formula \eqref{jsggsx1}. To illustrate this point, later we extract the main steps to prove the theorem from the reference \cite{Siegel}.

\section{Proof of main theorems}

\subsection{Proof of Theorem \ref{xdl2}}
\begin{proof}
Let $R$ be the positively oriented rectangular contour with $D(T)$. By the argument principle we get

\begin{align}\label{hljfa}
N(T)=\frac{1}{2\pi}\triangle_R\arg{\xi(s)}.
\end{align}
$\triangle_R\arg{\xi(s)}$ counts the changes in the argument of $\xi(s)$ along the contour $R$. We divide $R$ into three sub-contours. Let $L1$ be the horizontal line from
$-1$ to $2$. Let $L2$ be the sub-contour from $2$ to $2+iT$ and then to $1/2+iT$. Finally, let $L3$ be the sub-contour from $1/2+iT$ to $-1+iT$ and then $-1$. Accordingly\cite{ref19},

$$\triangle_R\arg{\xi(s)}=\triangle_{L_1}\arg{\xi(s)}+\triangle_{L_2}\arg{\xi(s)}+\triangle_{L_3}\arg{\xi(s)}$$
By $\xi(s)=\xi(1-s)$ we have $\xi(s)=\overline{\xi(1-\overline{s})}$. It shows that
$$\triangle_{L_2}\arg{\xi(s)}=\triangle_{L_3}\arg{\xi(s)}$$
Along the $L_1$ since $\xi(s)$ is a real function and $\triangle_{L_1}\arg{\xi(s)}=0$. Thus
\begin{align*}
N(T)&=\frac{1}{\pi}\triangle_{L_2}\arg{\xi(s)}.
\end{align*}
It is known that
\begin{align*}
\triangle_{L_2}\arg{\xi(s)}&=\triangle_{L_2}\arg{[\frac{s}{2}(s-1)\pi^{-s/2}\Gamma(\frac{s}{2})\zeta(s)]}.
\end{align*}
It follows that
\begin{align*}
\triangle_{L_2}\arg{[\frac{s}{2}(s-1)]}&=\arg{(1/2+iT)}+\arg{(-1/2+iT)}\\
&=\pi.
\end{align*}
and
\begin{align*}
\triangle_{L_2}\arg{[\pi^{-s/2}\Gamma(\frac{s}{2})\zeta(s)]}&=\triangle_{L_2}\arg{[\pi^{-s/2}\Gamma(\frac{s}{2})]}+\triangle_{L_2}\arg{\zeta(s)}\\
&=\arg{[\pi^{-\frac{1/2+iT}{2}}\Gamma(\frac{1/2+iT}{2})]}+\arg{\zeta(1/2+iT)}\\
&=\theta(T)+\arg{\zeta(1/2+iT)}.
\end{align*}
By the two above and \eqref{hljfa}, we obtain
\begin{align*}
N(T)&=\frac{1}{\pi}\triangle_{L_2}\arg{\xi(s)}\\
&=\frac{\theta(T)+\arg\zeta(1/2+iT)}{\pi}+1.
\end{align*}
This completes the proof.
\end{proof}

\subsection{Proof of Theorem  \ref{xdl3}}

\begin{proof}
From this formula \eqref{xgs1}, we get
\begin{align*}
\pi^{-\frac{s}{2}}\Gamma(\frac{s}{2})\zeta(s)&=\underline{-\frac{1}{s}+\int_1^{\infty}\psi(x)(x)x^{\frac{s}{2}-1}dx}\quad\underline{-\frac{1}{1-s}+\int_1^{\infty}\psi(x)(x)x^{-\frac{s+1}{2}}dx}.
\end{align*}
Let
$$\varphi(s)=-\frac{1}{s}+\int_1^{\infty}\psi(x)x^{\frac{s}{2}-1}dx,$$
then
\begin{align}\label{unknow}
\pi^{-\frac{s}{2}}\Gamma(\frac{s}{2})\zeta(s)=\varphi(s)+\varphi(1-s).
\end{align}
When $s=1/2+it$, $\pi^{-\frac{s}{2}}\Gamma(\frac{s}{2})\zeta(s)$ is a real function, and
\begin{align}\label{llld}
{\pi^{-\frac{s}{2}}\Gamma(\frac{s}{2})\zeta(s)}&=2\Re\left(\varphi(s)\right)\nonumber\\
&={2r(t)\cos[\omega(t)]},
\end{align}
$\mathrm{where} \;r(t)=|\varphi(s)|, \omega(t)=\arg\varphi(s).$

We know taht when $s=1/2+iT$, the included angle between $\varphi(s)$ and $\pi^{-\frac{s}{2}}\Gamma(\frac{s}{2})\zeta(s)$ meets
\begin{align*}
\left|\arg\varphi(s)- \arg\pi^{-\frac{s}{2}}\Gamma(\frac{s}{2})\zeta(s)\right|<{\frac{\pi}{2}},
\end{align*}
\begin{align}\label{fsx}
-\frac{\pi}{2}<\arg\varphi(s)-\arg\pi^{-\frac{s}{2}}\Gamma(\frac{s}{2})-\arg\zeta(s)<{\frac{\pi}{2}}.
\end{align}
According to \eqref{llld}, there are two types of zeros for $\zeta(1/2+it)$ function. One is
$\cos[\omega(t)]=0$, and the number of zeros is represented by
\begin{align*}
N_{\cos}(T)&=\#\{t:\cos[\omega(t)]=0,0\leq{t}\leq{T}\}.
\end{align*}
The other is $r(t)=0$, and the number of zeros is recorded as
\begin{align*}
N_{r}(T)&=\#\{t:r(t)=0\quad{ and}\quad \cos[\omega(t)]\neq{0},0\leq{t}\leq{T}\}.
\end{align*}
In this way, the number of zeros of $\zeta(1/2+it)$ function on $L(T)$
\begin{align}\label{N0Tsum}
N_0(T)=N_{\cos}(T)+N_r(T).
\end{align}

$\cos[\omega(t)]$ is a composite function. Let's suppose that $\omega(t)$ is a strictly increasing function. It is not difficult to prove that when $\omega(t)$ is strictly increasing, the number of roots of $\cos[\omega(t)]=0$ is the least in the interval $t\in{[0,T]}$, or in the interval $\omega(t)\in{[\omega(0),\omega(T)]}$. Next, we calculate the number of zeros of $\cos[\omega(t)]=0$ in the interval $\omega(t)\in{[\omega(0),\omega(T)]}$.

It is known that $t_1\approx{14.13}$ is the first zero of $\zeta(s)$ function. According to the theorem \ref{xdl2}, when $T<t_1$, we have $N(T)=0$. Take $T=0<t_1$, and substitute it into the formula \eqref{NTF1} to get
$$N(T)=\frac{\theta(0)+\arg{\zeta(\frac{1}{2}+0i)}}{\pi}+1=0,$$
\begin{align}\label{tmpf1}
\theta(0)+\arg{\zeta(\frac{1}{2}+0i)}=-\pi.
\end{align}
Substitute $T=0$ into the inequality \eqref{fsx}, we have
$$-\frac{\pi}{2}<\arg\varphi(\frac{1}{2}+0i)- \theta(0)-\arg{\zeta(\frac{1}{2}+0i)}<\frac{\pi}{2},$$
Substitute \eqref{tmpf1} into the above, then
$$-\frac{3\pi}{2}<\arg\varphi(\frac{1}{2}+0i)<-\frac{\pi}{2}.$$
By $\omega(t)=\arg\varphi(\frac{1}{2}+it)$, we get
\begin{align}\label{omegabds1}
-\frac{3\pi}{2}<\omega(0)<-\frac{\pi}{2}.
\end{align}
It is known that $t_2\approx{21.02}$ and $t_3\approx{25.01}$ are the second and third zeros of $\zeta(s)$ function, respectively. We take $T_1\in{(t_2, t_3)}$, and substitute it into formula \eqref{NTF1}, then
$$N(T_1)=\frac{\theta(T_1)+\arg{\zeta(\frac{1}{2}+iT_1)}}{\pi}+1=2,$$
\begin{align}\label{tmpf2}
\theta(T_1)+\arg{\zeta(\frac{1}{2}+iT_1)}&=\pi.
\end{align}
Putting $T_1$ into the inequality \eqref{fsx}, we have
$$-\frac{\pi}{2}<\arg\varphi(\frac{1}{2}+iT_1)- \theta(T_1)-\arg{\zeta(\frac{1}{2}+iT_1)}<\frac{\pi}{2},$$
Put \eqref{tmpf2} is substituted into the above, we obtain
$$\frac{3\pi}{2}<\arg\varphi(\frac{1}{2}+iT_1)<\frac{3\pi}{2},$$
namely
\begin{align}\label{omegabds2}
\frac{\pi}{2}<\omega(T_1)<\frac{3\pi}{2}.
\end{align}
Consider Formula \eqref{omegabds1} and Formula \eqref{omegabds2}. According to the intermediate value theorem, there must be $T_0\in[0,T_1]$ satisfying $\omega(T_0)=0$. Taking $T_0$ as the dividing point, divide the interval $[0, T]$ of the independent variable $t$ into
$$[0,T]=[0,T_0)\cup[T_0,T],$$
as well as
$$[\omega(0), \omega(T)]=[\omega(0),0)\cup[0,\omega(T)].$$
\begin{itemize}
\item\textbf{Calculate the number of zeros of $\cos[\omega(t)]$ in the interval
$\omega(t)\in[\omega (0),0)$.}
\end{itemize}
In view of Formula \eqref{omegabds1}. Because $-\frac{\pi}{2}\in[\omega(0),0)$, there must be $\exists{t_c\in[0,T_0)}$ satisfying $$\omega(t_c)=-\frac{\pi}{2},$$
and
$$\cos[\omega(t_c)]=\cos[-\frac{\pi}{2}]=0.$$
That is, $\zeta(s)$ function has at least a zero point of $s=1/2+it_c$ in the interval
$[\omega(0),0)$.
\begin{itemize}
\item\textbf{Calculate the number of zeros of $\cos[\omega(t)]$ in the interval
$\omega(t)\in{[0,\omega(T)]}$.}
\end{itemize}

According to $\cos(k\pi-\frac{\pi}{2})=0$ and $\cos[\omega(t)]=0$, we get the formula to calculate the number of zeros of $\cos[\omega(t)]$ in the interval $\omega(t)\in{[0, \omega(T)]}$
\begin{align}\label{kcode}
k&=\left\lfloor\frac{\omega(T)}{\pi}+\frac{1}{2}\right\rfloor
\end{align}
Let $s=1/2+iT$, from \eqref{fsx}, we get
$$\frac{\theta(T)+\arg\zeta(1/2+iT)}{\pi}<\frac{\arg\varphi(1/2+iT)}{\pi}+\frac{1}{2}<\frac{\theta(T)+\arg\zeta(1/2+iT)}{\pi}+1.$$
Substituting $\arg\varphi(1/2+iT)=\omega(T)$ into the above, we have
$$\frac{\theta(T)+\arg\zeta(1/2+iT)}{\pi}<\frac{\omega(T)}{\pi}+\frac{1}{2}<\frac{\theta(T)+\arg\zeta(1/2+iT)}{\pi}+1,$$
and
$$\left\lfloor\frac{\omega(T)}{\pi}+\frac{1}{2}\right\rfloor=\frac{\theta(T)+\arg{\zeta(1/2+iT)}}{\pi}.$$
Putting the above into formula \eqref{kcode}, we get
$$k=\frac{\theta(T)+\arg{\zeta(1/2+iT)}}{\pi}$$
Therefore, the total number of zeros of $\cos[\omega(t)]$ in the interval $[\omega(0), 0)\cup[0, \omega(T)]$ is
$$N_{\cos}(T)={k+1}=\frac{\theta(T)+\arg{\zeta(1/2+iT)}}{\pi}+1.$$
Putting $N_{cos}(T)$ into formula \eqref{N0Tsum}, we obtain
\begin{align*}
N_0(T)=\frac{\theta(T)+\arg{\zeta(1/2+iT)}}{\pi}+1+N_r(T)
\end{align*}
Putting formula \eqref{NTF1} into the above, we have
$$N_0(T)=N(T)+N_r(T).$$
From $N(T)\geq{N_0(T)}$ and $N_r(T)\geq{0}$, we get $N_{r}(T)=0$. So
\begin{align*}
N_0(T)&=N(T)\\
&=\frac{\theta(T)+\arg{\zeta(1/2+iT)}}{\pi}+1.
\end{align*}
\end{proof}

\subsection{Proof of Theorem \ref{xdl4}}
\begin{proof}
By(Detail to see \cite{Siegel})
  $$\Phi(u)=\int\frac{e^{\pi{i}{x^2}+2\pi{iux}}}{e^{\pi{ix}}-e^{-\pi{ix}}}dx$$
one gets
\begin{align}\label{zb11}
\int_{0\nwarrow{1}}\frac{e^{-\pi{i}x^2+2\pi{iux}}}{e^{\pi{ix}}-e^{-\pi{ix}}}dx&=\frac{1}{1-e^{-2\pi{iu}}}-\frac{e^{\pi{iu^2}}}{e^{\pi{iu}}-e^{-\pi{iu}}}
\end{align}
Multiply both sides of the formula \eqref{zb11} by ${u^{-s}du}$, and integrate along the ray from $u=0$ to $u=i^{1/2}\infty$, $\underline{u^{-s}\textbf{ is defined on the slit plane}}$(excluding 0 and $-\infty$)
$$
\int_{0}^{\overline{\varepsilon}\infty}u^{-s}\int_{0\nwarrow{1}}\frac{e^{-\pi{i}x^2+2\pi{iux}}}{e^{\pi{ix}}-e^{-\pi{ix}}}dxdu=
\int_{0}^{\overline{\varepsilon}\infty}u^{-s}\frac{1}{1-e^{-2\pi{iu}}}du-\int_{0}^{\overline{\varepsilon}\infty}u^{-s}\frac{e^{\pi{iu^2}}}{e^{\pi{iu}}-e^{-\pi{iu}}}du
$$
where $\overline{\varepsilon}=i^{1/2}$. From the above, Riemann obtained
\begin{align}\label{zb12}
\pi^{-\frac{s}{2}}\Gamma(\frac{s}{2})\zeta(s)&=
\pi^{-\frac{s}{2}}\Gamma(\frac{s}{2})\int_{0\swarrow{1}}\frac{x^{-s}e^{\pi{ix^2}}}{e^{\pi{ix}}-e^{-\pi{ix}}}dx+
\pi^{-\frac{1-s}{2}}\Gamma(\frac{1-s}{2})\int_{0\searrow{1}}\frac{e^{-\pi{ix^2}}x^{s-1}}{e^{\pi{ix}}-e^{-\pi{ix}}}dx
\end{align}
One now sets
\begin{align}\label{syy40}
f(s)&=\int_{0\swarrow{1}}\frac{x^{-s}e^{\pi{ix^2}}}{e^{\pi{ix}}-e^{-\pi{ix}}}dx,
\end{align}
and put it into the \eqref{zb12},one gets
\begin{align}\label{zb12a}
{{\pi^{-\frac{s}{2}}\Gamma(\frac{s}{2})\zeta(s)=\pi^{-\frac{s}{2}}\Gamma(\frac{s}{2})f(s)+\pi^{-\frac{1-s}{2}}\Gamma(\frac{1-s}{2})f(1-s)}}
\end{align}
We now set
$$\varphi(s)=\pi^{-\frac{s}{2}}\Gamma(\frac{s}{2})f(s),$$
put it into the formula \eqref{zb12a}, then
\begin{align}\label{zb13}
\pi^{-\frac{s}{2}}\Gamma(\frac{s}{2})\zeta(s)&=\varphi(s)+\varphi(1-s)
\end{align}
When $s=\frac{1}{2}+it$, Riemann gives
\begin{align*}
\pi^{-\frac{s}{2}}\Gamma(\frac{s}{2})\zeta(s)=2\Re{[\varphi(s)]}\quad,\qquad{s=\frac{1}{2}+it}
\end{align*}
\end{proof}
\section{Acknowledgements}
I would like to express acknowledgement here to the Delaware State University Mr. Yi Ling for his help.


\small
\begin{thebibliography}{9999}
\bibitem{riemann} Bernhard Riemann. \"{U}ber die Anzahl der Primzahlen unter einer gegebenen Grosse.
Monatsberichte der Berliner Akademie, 1859.
\bibitem{Siegel} C.L. Siegel, \"{U}ber Riemanns Nachlass zur analytischen Zahlentheorie,Quellen und Studien zur Geschichte der Mathematik, Astronomie und Physik 2 (1932), also Gesammelte Abhandlungen, Bd. I, Springer-Verlag, Berlin.
\bibitem{SiegelX} https://arxiv.org/abs/1810.05198v1
\bibitem{Conrey} J. B. Conrey, More than two fifths of the zeros of the Riemann zeta function are on the critical line, J.Reine Angew. Math.
\bibitem{ref3} Juan Arias de Reyna and Jan van de Lune,A first encounter with the Riemann Hypothesis and its numerical verification, La Gaceta de la RSME, Num.
\bibitem{levins} N. Levinson, More than One Third of Zeros of Riemann's Zeta-Function are on G = l/2,ADVANCES IN MATHEMATICS.
\bibitem{ref32} Norman Levinson, Remarks on a Formula of Riemann for his Zeta-function, Journal of mathematical analysys and applications.
\bibitem{Titch} E.C.TITCHMARSH.The theory of the Riemann Zeta Function. Clarendon press oxford.
\bibitem{xref1} Harold M. Edwards. Riemann's Zeta Function. Dover, 1974.
\bibitem{ref19} Peter Borwein, Stephen Choi, Brendan Rooney and Andrea Weirathmueller. The Riemann Hypothesis. July 12, 2006
\bibitem{ref59} Hans Rademacher.Topics in Analytic Number Theory. Springer-Verlag Berlin, 1973
\bibitem{Selberg1} Selberg, A., On the zeros of Riemann's zeta-function. Skr. Norske Vid.-Akad. Oslo No.10(1942).
\end{thebibliography}
\end{document}